\documentclass[fleqn]{mat01}
\usepackage{times,mathtimy,amssymb}
\begin{document}

\setcounter{page}{213} \firstpage{213}

\newtheorem{theoreee}{Theorem}
\renewcommand\thetheoreee{\Alph{theoreee}}
\newtheorem{theore}[theoreee]{\bf Theorem}

\newtheorem{lem}{Lemma}

\def\defini{\trivlist \item[\hskip \labelsep{DEFINITION}]}
\def\theor{\trivlist \item[\hskip \labelsep{\bf Theorem.}]}
\def\pot{\trivlist \item[\hskip \labelsep{\it Proof of Theorem.}]}
\def\remar{\trivlist \item[\hskip \labelsep{\it Remark.}]}

\title{Law of iterated logarithm for NA sequences with non-identical
distributions}

\markboth{Guang-hui Cai and Hang Wu}{Law of iterated logarithm for
NA sequences}

\author{GUANG-HUI CAI and HANG WU$^{*}$}

\address{Department of Mathematics; $^{*}$Department of Statistics, Zhejiang
Gongshang University, Hangzhou 310035, People's Republic China\\
\noindent E-mail: cghzju@163.com; whwhzju@163.com}

\volume{117}

\mon{May}

\parts{2}

\pubyear{2007}

\Date{MS received 5 June 2005; revised 5 January 2007}

\begin{abstract}
Based on a law of the iterated logarithm for independent random
variables sequences, an iterated logarithm theorem for NA
sequences with non-identical distributions is obtained. The proof
is based on a Kolmogrov-type exponential inequality.
\end{abstract}


\maketitle

\section{Introduction}

There are many results about law of iterated logarithm for
independent random variables sequences. The following theorem is
one of them.

\begin{theore}[({Theorem 7.2 of \cite{6}})]
Let $\{X_{n},n\geq 1\}$ be independent random variables
sequences{\rm ,} $EX_{n}=0, \sigma_{n}^{2}=EX_{n}^{2}<\infty,
B_{n}'=\sum_{k=1}^{n}\sigma_{k}^{2}, S_{n}=\sum_{k=1}^{n}X_{k},$
$\triangle_{n}=\sup_{x}$ $|P(S_{n}<x\sqrt{B_{n}'})-\Phi(x)|,$
where $\Phi(x)$ is a standard normal distribution function.
If\vspace{-.5pc}
\begin{enumerate}
\renewcommand\labelenumi{\rm (\roman{enumi})}
\leftskip .4pc
\item $B_{n}'\longrightarrow \infty,$ when
$n\longrightarrow\infty,$

\item $B_{n+1}'/B_{n}'\longrightarrow 1,$ when
$n\longrightarrow\infty,$

\item $\triangle_{n}=O[(\log B_{n}')^{-1-\delta}],\delta>0,$
\end{enumerate}\vspace{-1pc}
hold{\rm ,} then
\begin{equation*}
\hskip -1.25pc
\limsup_{n\longrightarrow\infty}\frac{S_{n}}{(2B_{n}' \log\!\log
B_{n}')^{{1}/{2}}}=1 \ \hbox{\rm a.s.}
\end{equation*}
\end{theore}

As for negatively associated (NA) random variables, Joag \cite{2}
gave the following definition.

\begin{defini}{\cite{2}}$\left.\right.$\vspace{.5pc}

\noindent A finite family of random variables $\{X_{i},1\leq i\leq
n\} $ is said to be negatively associated (NA) if for every pair
of disjoint subsets $T_{1}$ and $T_{2}$ of $\{1,2,\dots,n\}$, we
have
\begin{equation*}
\hbox{Cov}(f_{1}(X_{i},i\in T_{1}),f_{2}(X_{j},j\in T_{2}))\leq 0,
\end{equation*}
whenever $f_{1}$ and $f_{2}$ are coordinatewise increasing and the
covariance exists. An infinite family is negatively associated if
every finite subfamily is negatively associated.
\end{defini}

Recently, some authors focused on the problem of limiting behavior
of partial sums of NA sequences. Su {\it et~al} \cite{9} derived
some moment inequalities of partial sums and a weak convergence
for a strong stationary NA sequence. Lin \cite{5} set up an
invariance principal for NA sequences. Su and Qin \cite{10} also
studied some limiting results for NA sequences. More recently,
Liang \cite{3,4} considered some complete convergence for weighted
sums of NA sequences. Those results, especially some moment
inequality by Huang and Xu \cite{1}, Shao \cite{7} and Yang
\cite{11}, undoubtedly propose important theory guide in further
apply for the NA sequence. Shao \cite{8} obtained a law of
iterated logarithm for NA sequences with identical distributions.
Zhang \cite{12} got a law of iterated logarithm for NA vector with
identical  distributions.

Based on a law of the iterated logarithm for independent random
variables sequences, the main purpose of this paper is to
establish an iterated logarithm theorem for NA sequences with
non-identical distributions. The proof is based on a Kolmogrov
type exponential inequality. Hence the following theorem.

\begin{theor}{\it
Let $\{X_{n},n\geq 1\}$ be NA random variables sequences{\rm ,}
$EX_{n}=0,$ $\sup_{n\geq 1}$ $EX_{n}^{2}(\log| X_{n}|)^{1+\delta}
<\infty,$ for some $\delta>0.$ Let $S_{n}=\sum_{k=1}^{n}X_{k},$
$B_{n}=\hbox{\rm Var}\, S_{n}>0, B_{n}'=\sum_{k=1}^{n}EX_{k}^{2},
$ $\triangle_{n}=\sup_{x}|P(S_{n}<x\sqrt{B_{n}})-\Phi(x)|,$ where
$\Phi(x)$  is a standard normal distribution function.
If\vspace{-.5pc}
\begin{enumerate}
\renewcommand\labelenumi{\rm (\roman{enumi})}
\leftskip .4pc
\item $B_{n}=O(n),$

\item $B_{n+1}/B_{n}\longrightarrow 1,$ when
$n\longrightarrow\infty,$

\item $\triangle_{n}=O[(\log B_{n})^{-1}],$

\item $B_{n}/B_{n}'\longrightarrow 1,$ when
$n\longrightarrow\infty,$
\end{enumerate}\vspace{-1pc}
hold{\rm ,} then
\begin{equation}
\hskip -1.25pc
\limsup_{n\longrightarrow\infty}\frac{S_{n}}{(2B_{n}\log\!\log
B_{n})^{{1}/{2}}}=1 \ \hbox{a.s.}
\end{equation}\vspace{-2pc}}
\end{theor}

\section{The proof of theorem}

Throughout this paper, $C$ will represent a positive constant
though its value may change from one appearance to the next, and
$a_{n}=O(b_{n})$ will mean $a_{n}\leq C b_{n}$, and $a_{n}\ll
b_{n}$ will mean $a_{n}=O(b_{n}).$

In order to prove our results, we need the following lemma.

\begin{lem}\hskip -.4pc{\rm \cite{7}.} \ \
Let $\{X_{n},n\geq 1\}$  be  NA  random variables sequences{\rm ,}
$EX_{n}=0, EX_{n}^{2}<\infty, B_{n}=\sum_{k=1}^{n}EX_{k}^{2},
S_{n}=\sum_{k=1}^{n}X_{k}$. Then for all  $x>0, a>0$  and
$0<\alpha<1,$ we have
\begin{align*}
&P\Big(\max_{1\leq i\leq n}|S_{i}|\geq x\Big)\\[.5pc]
&\quad\, \leq 2P\Big(\max_{1\leq i\leq n}|X_{i}|>a
\Big)+\frac{2}{1-\alpha}\exp\left\{-\frac{x^{2}\alpha}{2(ax+
B_{n})}\right\}, \quad n\geq 1.
\end{align*}
\end{lem}

\begin{pot}{\rm For all
$i\geq 1$, let $X_{i}^{(n)}=X_{i}I(|X_{i}|\leq a_{i})+
a_{i}I(X_{i}>a_{i})-a_{i}I(X_{i}<-a_{i})$,
$\forall\varepsilon\in\big(0,\frac{1}{20}\big),$ where
$a_{i}=\varepsilon\big(\frac{B_{i}}{\log\!\log
B_{i}}\big)^{\frac{1}{2}}.$ Let
$S_{n,1}=\sum_{i=1}^{n}(X_{i}^{(n)}-EX_{i}^{(n)})$,
$S_{n,2}=\sum_{i=1}^{n}[X_{i}-X_{i}^{(n)}-E(X_{i}-X_{i}^{(n)})]$.
By $\sup_{n\geq 1}EX_{n}^{2}(\log|X_{n}|)^{1+\delta}<\infty$ and
$B_{n}=O(n),$ we have
\begin{align*}
&\sum_{i=1}^{\infty}E|X_{i}|I(|X_{i}|> a_{i})/(B_{i}\log\!\log
B_{i})^{{1}/{2}}\nonumber\\[.3pc]
&\quad\,\leq \sum_{i=1}^{\infty}EX_{i}^{2}(\log|X_{i}|)^{1+\delta}I(|X_{i}|> a_{i}) \frac{a_{i}^{-1} (\log a_{i})^{-1-\delta}}{(B_{i}\log\!\log B_{i})^{\frac{1}{2}}}  \nonumber\\[.3pc]
&\quad\,\ll \sum_{i=1}^{\infty} \frac{1}{(B_{i}\log\!\log B_{i})^{{1}/{2}}\cdot \left(\frac{B_{i}}{\log\!\log B_{i}}\right)^{{1}/{2}}}\cdot \left[\log\left(\frac{B_{i}}{\log\!\log B_{i}}\right)^{{1}/{2}} \right]^{-1-\delta} \nonumber\\[.3pc]
&\quad\,\ll \sum_{i=1}^{\infty} \frac{1}{i (\log
i)^{1+\frac{\delta}{2}} }<\infty.
\end{align*}}
\end{pot}
By Kronecker lemma, we have
\begin{equation*}
\frac{\sum\limits_{i=1}^{\infty}E|X_{i}|I(|X_{i}|>
a_{i})+\sum\limits_{i=1}^{\infty}|X_{i}|I(|X_{i}|>
a_{i})}{(B_{n}\log\!\log B_{n})^{{1}/{2}}}\longrightarrow 0 \
\hbox{a.s.}
\end{equation*}
Then we have
\begin{equation}
\frac{S_{n,2}}{(B_{n}\log\!\log B_{n})^{{1}/{2}}}\longrightarrow 0
\ \hbox{a.s.}
\end{equation}
By Lemma~1 and $\alpha=1-\varepsilon, a=2a_{n}$, we have
\begin{align*}
&P\Big(\max_{1\leq i\leq n}|S_{i,1}|>(1+4\varepsilon)(2B_{n}\log\!\log B_{n})^{{1}/{2}}\Big)\nonumber\\[.3pc]
&\quad\, \leq
\frac{2}{\varepsilon}\exp\left(-\frac{(1+4\varepsilon)^{2}(2B_{n}\log\!\log
B_{n})(1-\varepsilon)}{4a_{n}(1+4\varepsilon)(2B_{n}\log\!\log
B_{n})^{{1}/{2}}+2B_{n}'}\right).
\end{align*}
So we have
\begin{align}
&\sum_{n=1}^{\infty}\frac{1}{n}P\Big(\max_{1\leq i\leq n}|S_{i,1}|>(1+4\varepsilon)(2B_{n}\log\!\log B_{n})^{{1}/{2}}\Big) \nonumber\\[.3pc]
&\quad\,\leq\sum_{n=1}^{\infty}\frac{1}{n}\cdot\frac{2}{\varepsilon}\exp\left(-\frac{(1+4\varepsilon)^{2}(2B_{n}\log\!\log B_{n})(1-\varepsilon)}{4a_{n}(1+4\varepsilon)(2B_{n}\log\!\log B_{n})^{{1}/{2}}+2B_{n}'}\right)\nonumber\\[.3pc]
&\quad\,\ll\sum_{n=1}^{\infty}\frac{1}{n}\cdot\frac{2}{\varepsilon}\exp\left(-\frac{(1+4\varepsilon)^{2}(1-\varepsilon) \log\!\log B_{n}}{2\sqrt{2}\varepsilon (1+4\varepsilon)+1+\varepsilon}\right) \nonumber\\[.3pc]
&\quad\,\leq\sum_{n=1}^{\infty}\frac{1}{n}\cdot\frac{2}{\varepsilon}\exp(-(1+\varepsilon) \log\!\log B_{n})\nonumber
\end{align}\pagebreak
\begin{align}
&\quad\,\ll\sum_{n=1}^{\infty}\frac{1}{n}\exp(-(1+\varepsilon) \log\!\log n)  \nonumber\\[.3pc]
&\quad\,=\sum_{n=1}^{\infty}\frac{1}{n(\log
n)^{1+\varepsilon}}<\infty.
\end{align}
By (3), for all $\varepsilon>0$, we have
\begin{align*}
\sum_{n=1}^{\infty}\frac{1}{n}P\Big(\max_{1\leq i\leq
n}|S_{i,1}|>(1+4\varepsilon)(2B_{n}\log\!\log
B_{n})^{{1}/{2}}\Big) <\infty.
\end{align*}
Then, we have
\begin{align*}
\hskip
-4pc\sum\limits_{k=0}^{\infty}\sum\limits_{n=2^{k}}^{2^{k+1}-1}
(2^{k+1}-1)^{-1}P\Big( \max\limits_{1\leq j\leq
2^{k}}|S_{j,1}|>(1+4\varepsilon)(2B_{2^{k+1}-1}\log\!\log
B_{2^{k+1}-1})^{{1}/{2}}\Big)<\infty.
\end{align*}
Then
\begin{equation*}
\sum\limits_{k=0}^{\infty} P\Big( \max\limits_{1\leq j\leq
2^{k}}|S_{j,1}|>(1+4\varepsilon)(2B_{2^{k+1}-1}\log\!\log
B_{2^{k+1}-1})^{{1}/{2}}\Big)<\infty.
\end{equation*}
So, we have
\begin{align*}
&\frac{\max\limits_{1\leq j\leq
2^{k}}|S_{j,1}|}{(2B_{2^{k+1}-1}\log\!\log
B_{2^{k+1}-1})^{{1}/{2}}}\leq 1+4\varepsilon \ \hbox{a.s.}
\end{align*}
Notice that for any positive integer $n$ there exists a
non-negative integer $k_{0}$, such that $2^{k_{0}}\leq n<
2^{k_{0}+1}$. Then, we have
\begin{align*}
\frac{|S_{n,1}|}{(2B_{n}\log\!\log B_{n})^{{1}/{2}}} \leq \frac{
\max\limits_{1\leq j\leq
2^{k_{0}+1}}|S_{j,1}|}{(2B_{2^{k_{0}+1}-1}\log\!\log
B_{2^{k_{0}+1}-1})^{{1}/{2}}} \leq 1+4\varepsilon \ \hbox{a.s.}
\end{align*}
So, we have
\begin{equation}
\limsup_{n\longrightarrow
\infty}\frac{|S_{n,1}|}{(2B_{n}\log\!\log B_{n})^{{1}/{2}}}\leq
1+4\varepsilon \ \hbox{a.s.}
\end {equation}
By $S_{n}=S_{n,1}+S_{n,2}$, (2) and (4),
\begin{equation}
\limsup_{n\longrightarrow \infty}\frac{S_{n}}{(2B_{n}\log\!\log
B_{n})^{{1}/{2}}}\leq 1+5\varepsilon \ \hbox{a.s.}
\end{equation}
Next, we prove the following:
\begin{equation}
\limsup_{n\longrightarrow \infty}\frac{S_{n}}{(2B_{n}\log\!\log
B_{n})^{{1}/{2}}}> 1-\varepsilon \ \hbox{a.s.}
\end{equation}
By Theorem~1(iii), we have
$\triangle_{n}=\sup_{x}|1-\Phi(x)-P(S_{n}> x\sqrt{B_{n}})|$. Then
\begin{align}
&P(S_{n}>(1-\varepsilon)(2B_{n}\log\!\log B_{n})^{{1}/{2}}) \nonumber\\[.3pc]
&\quad\,\geq 1-\phi((1-\varepsilon)(2\log\!\log B_{n})^{{1}/{2}})- \triangle_{n}\nonumber\\[.3pc]
&\quad\,\geq C\frac{1}{\sqrt{2\pi}(1-\varepsilon)(2\log\!\log B_{n})^{{1}/{2}}}\nonumber\\[.6pc]
&\qquad\, \times \exp \left(-\frac{(1-\varepsilon)^{2}(2\log\!\log B_{n})}{2}\right) -C(\log B_{n})^{-1}\nonumber\\[.3pc]
&\quad\,\geq \frac{C}{(\log
B_{n})^{(1-\varepsilon)^{2}}(\log\!\log B_{n})^{{1}/{2}}}.
\end{align}
For all $\tau>0,$ there exists non-decreasing positive integers
sequence $\{n_{k}\}$. We have $n_{k}\longrightarrow\infty$ when
$k\longrightarrow\infty,$ and
\begin{align}
B_{n_{k}-1}\leq (1+\tau)^{k}<B_{n_{k}}, \quad k=1,2,\dots.
\end{align}
Let $\psi(n_{k})=(2(B_{n_{k}}-B_{n_{k-1}})\log\!\log
(B_{n_{k}}-B_{n_{k-1}}))^{{1}/{2}},$ for all $r\in (0,1).$ By (7),
we have
\begin{align}
&\sum_{k=1}^{\infty}P((S_{n_{k}}-S_{n_{k-1}})>(1-r)\psi(n_{k})) \nonumber\\[.3pc]
&\quad\,
\geq\sum_{k=1}^{\infty}\frac{C}{(\log(B_{n_{k}}-B_{n_{k-1}}))^{(1-r)^{2}}(\log\!\log
(B_{n_{k}}-B_{n_{k-1}}))^{{1}/{2}}}\nonumber\\[.3pc]
&\quad\, \geq\sum_{k=1}^{\infty}\frac{C}{k^{(1-r)^{2}}(\log
k)^{{1}/{2}}}=\infty.
\end{align}
By the generalized Borel--Cantelli lemma and (9), we have
\begin{align}
P((S_{n_{k}}-S_{n_{k-1}})>(1-r)\psi(n_{k}) \ \hbox{i.o.})=1.
\end{align}
By (7), when $n$ large enough, we have
\begin{align}
|S_{n}|\leq 2(2B_{n}\log\!\log B_{n})^{{1}/{2}} \ \hbox{a.s.}
\end{align}
Let $\chi(n)=(2B_{n}\log\!\log B_{n})^{{1}/{2}}$, when
$k\longrightarrow\infty$,
\begin{align}
(1\!-\!r)\psi(n_{k})-2\chi(n_{k-1})\sim
[(1\!-\!r)\tau^{{1}/{2}}(1+\tau)^{-{1}/{2}}-2(1+
\tau)^{-{1}/{2}}]\chi(n_{k}).
\end{align}
Because $\varepsilon \in R^{+}$, we can choose $r$ and $\tau$.
Then
\begin{align}
(1-r)\tau^{{1}/{2}}(1+\tau)^{-{1}/{2}}-2(1+
\tau)^{-{1}/{2}}>1-\varepsilon.
\end{align}

$\left.\right.$\vspace{-1.5pc}

\pagebreak

\noindent By (12), (13), (11) and (10), we have
\begin{align*}
&P(S_{n_{k}}>(1-\varepsilon)\chi(n_{k}) \ \hbox{i.o.})\nonumber\\[.4pc]
&\quad\, \geq P(S_{n_{k}}>(1-r)\psi(n_{k})-2\chi(n_{k-1}) \ \hbox{i.o.})\nonumber\\[.4pc]
&\quad\, \geq P((S_{n_{k}}-S_{n_{k-1}})>(1-r)\psi(n_{k}) \
\hbox{i.o.})=1.
\end{align*}
Now we complete the proof of (6). By (5) and (6), (1) holds.

\begin{remar}
Theorem~1 generalizes the Kolmogrov type law of iterated logarithm
(see Theorem~7.2 of \cite{6}) to NA sequences.
\end{remar}

\section*{Acknowledgments}

The author would like to thank the referee for his/her valuable
comments. This paper is supported by Society Key Research
Foundation of Zhejiang Province University (Statistics of Zhejiang
Gongshang University) and key discipline of Zhejiang Province (key
discipline of  Statistics of Zhejiang Gongshang University).


\begin{thebibliography}{99}
\bibitem{1}Huang W T and Xu B, Some maximal inequalities
and complete convergences of negatively associated random
sequences, {\it Statist. Probab. Lett.} {\bf 57} (2002) 183--191

\bibitem{2}Joag D K and Proschan F, Negative associated of random
variables with application, {\it Ann. Statist.} {\bf 11} (1983)
286--295

\bibitem{3}Liang H Y and Su C, Complete convergence for weighted
sums of NA sequences, {\it Statist. Probab. Lett.} {\bf 45} (1999)
85--95

\bibitem{4}Liang H Y, Complete convergence for weighted sums of
negatively associated random variables, {\it Statist. Probab.
Lett.} {\bf 48} (2000) 317--325

\bibitem{5}Lin Z Y, Invariance principle for negatively
associated sequence, {\it Chinese Sci. Bull.} {\bf 42} (1997)
238--242 (in Chinese)

\bibitem{6}Petrov V V, Limit theorems of probability
theory~-- sequences of independent random variables (Oxford:
Oxford Science Publications) (1995)

\bibitem{7}Shao Q M, A comparison theorem on moment inequalities
between Negatively associated and independent random variables,
{\it J.~Theor. Probab.} {\bf 13} (2000) 343--356

\bibitem{8}Shao Q M and Su C, The law of the iterated logarithm
for negatively associated random variables, {\it Stochastic
Process Appl.} {\bf 83} (1999) 139--148

\bibitem{9}Su C, Zhao L C and Wang Y B Moment inequalities
and weak convergence for NA sequences, {\it Science in China} {\bf
A26} (1996) 1091--1099 (in Chinese)

\bibitem{10}Su C and Qin Y S, Limit theorems for negatively
associated sequences, {\it Chinese Sci. Bull.} {\bf 42} (1997)
243--246

\bibitem{11}Yang S C, Moment inequality of random variables
partial sums, {\it Science in Chinese} {\bf A30} (2000) 218--223

\bibitem{12}Zhang L X, Strassen's law of the iterated
logarithm for negatively associated random vectors, {\it
Stochastic Process Appl.} {\bf 95} (2001) 311--328
\end{thebibliography}
\end{document}